\documentclass[a4paper,10pt]{article}
\usepackage[english,activeacute]{babel}
\usepackage{amssymb,amsmath,amsthm,amsfonts,mathrsfs,scalefnt,graphicx,graphics,color,hyperref,}
\usepackage[T1]{fontenc}
\usepackage{url}  
\usepackage[latin1]{inputenc}
\theoremstyle{plain}
\newtheorem{thm}{Theorem}
\newtheorem{cor}[thm]{Corollary}
\theoremstyle{definition}

\theoremstyle{plain}
\newtheorem{lemme}[thm]{Lemma}
\newtheorem{prop}[thm]{Proposition}
\theoremstyle{remark}
\newtheorem{rmq}{Remark}
\newcommand{\Cb}  {{\mathbb C}}

\newcommand{\Rb}  {{\mathbb R}}

\DeclareMathOperator{\re}{Re}

\newcommand{\izi}{\int\limits_0^{+\infty}}

\newcommand{\id}[2]{\int\limits_{#1}^{#2}}

\newcommand{\ef}  {{\mathbf e}}
\newcommand{\tx}[1][x]{T_{[#1,+\infty[}}
\newcommand{\txh}[1][x]{T_{#1}^{h}}
\newcommand{\txz}[1][x]{T_{#1}^{0}}
\newcommand{\xtx}[1][x]{X_{T_{[#1,+\infty[}}}
\newcommand{\sel}[1][\lambda]{{S_{\ef_{#1}}}}
\newcommand{\se}{{S_{\ef}}}
\newcommand{\psip}[2]{\psi_{#1}^{+}(-{#2})}

\newcommand{\plx}{P_\lambda^{(x)}}
\newcommand{\elx}{E_\lambda^{(x)}}

\newcommand{\Ls} {{\mathcal L}}

\title{On the scaling property in fluctuation theory for stable Lévy processes}
\author{Fernando Cordero\protect\footnote{Laboratoire de Probabilités et Modèles Aléatoires, UMR 7599, UPMC, Case courrier 188, 4, Place Jussieu, 75252 PARIS Cedex 05, France. E-mail : \protect\url{cordiery@gmail.com}}}

\begin{document}

\maketitle

\begin{abstract}
We find an expression for the joint Laplace transform of the law of $(\tx,\xtx)$ for a Lévy process $X$, where 
$\tx$ is the first \mbox{hitting} time of $[x,+\infty[$ by $X$.
When $X$ is an $\alpha$-stable Lévy process, with $1<\alpha<2$, 
we show how to recover from this formula the law of $\xtx$; this result was already obtained by D. Ray,
in the symmetric case and by N. Bingham, in the case when $X$ is non spectrally negative. Then,  
we study the behaviour of the time of first passage $\tx$ conditioned to $\{\xtx -x \leq h\}$ when 
$h$ tends to $0$. This study brings forward an asymptotic variable $\txz$, which seems to be related to the absolute \mbox{continuity} of
the law of the supremum of $X$.\\
\end{abstract}


\section{Introduction}
\noindent Let $X=(X_t, t\geq 0)$ denote a Lévy process with characteristic exponent $\psi$.\\

\noindent For $t\geq 0$ we set :

 $$S_t=\sup\limits_{0\leq s\leq t}X_s\quad\textrm{ and }\quad I_t=\inf\limits_{0\leq s\leq t}X_s,$$

\noindent and for $x>0$ :
 $$T_{[x,+\infty[}=\inf\{s>0: X_s\geq x\}\textrm{ and } K_x=X_{T_{[x,+\infty[}}-x.$$\\

\noindent We denote by $\bf{e}_\gamma$ an exponential variable with parameter $\gamma$, independent of $X$.\\

\noindent Our interest focusses on the joint law of $(T_{[x,+\infty[},X_{T_{[x,+\infty[}})$ (see \cite{gu}) and on questions concerned with
 the absolute continuity of the law of $S_t$. \\ 

\noindent In section 2 we present some known results, which we will be needed in the following sections. 
The first result (Proposition \ref{proppero}) is a famous formula obtained by Pecherskii and Rogozin which concerns
the joint Laplace transform of $(\tx[\ef_\gamma],\xtx[\ef_\gamma])$. The second result 
(Proposition \ref{propber}) deals with the asymptotic behaviour of the quantity $P(S_1<x)$ in the stable case 
when $x$ tends to $0$.\\

\noindent In section 3 we obtain an expression of the joint Laplace transform of the pair
$(\tx,\xtx)$. First, we remark that the Laplace transform of this quantity can be expressed
by means of the joint Laplace transform of $(\tx[\ef_\gamma],\xtx[\ef_\gamma])$, 
for which we dispose of a formula (Proposition \ref{proppero}). In this way we obtain the result by inverting 
a Laplace transform.\\

\noindent Then, we pass to the $\alpha$-stable case, with $1<\alpha<2$. By means of 
the scaling property and of Proposition \ref{propber}, we prove some asymptotic properties associated to quantities introduced in the 
previous paragraph. After this, we recover the law of $\xtx$, as a corollary of these properties. This law is already known. 
In \cite{Ray} Ray gives an expression of its density in the symmetric stable case.
In \cite{Bin} Bingham generalises this result for non-spectrally negative stable processes. However, the interest of this part of the work 
is to point out the role of the scaling property in the obtention of this result. \\

\noindent Finally, we focus on the asymptotic behaviour of the random variable $\tx$ conditionned upon 
$\{K_x\leq h\}$ when $h$ tends to $0+$. We show that a convergence in law holds towards a random variable which we denote $\txz$.\\

\section{Preliminaries}
\noindent First, we introduce some notations that will be useful afterwards. For $q>0$, $\re{\lambda}\leq 0$ and $\re{\mu}\geq 0$, we define :
\begin{equation}
\psi_{q}^+(\lambda)=\exp\Bigg(-\int\limits_0^{+\infty} (e^{\lambda x}-1)d_x\Big( \int\limits_0^{+\infty} u^{-1}e^{-q u}P(X_u>x)du\Big)\Bigg),
\end{equation}
\noindent and
\begin{equation}
\psi_{q}^-(\mu)=\exp\Bigg(\int\limits_{-\infty}^0 (e^{\mu x}-1)d_x\Big( \int\limits_0^{+\infty} u^{-1}e^{-q u}P(X_u<x)du\Big)\Bigg).
\end{equation}
\noindent We use the following result due to Pecherskii and Rogozin :
\begin{prop}[\cite{pero}, p.420]\label{proppero}
 For all $\gamma,\lambda,\mu>0$, we have :\\
\begin{equation}\label{eftk1}
E[\exp(-\lambda T_{[\bf{e}_\gamma,+\infty[}-\mu K_{\bf{e}_\gamma})]=
\frac{\gamma}{\gamma-\mu}\bigg( 1-\frac{\psi_{{\lambda}}^+(-\gamma)}{\psi_{{\lambda}}^+(-\mu)}\bigg).
\end{equation}
\,
\end{prop}
\,

\begin{rmq}
It follows from the definition of $\psi_{{\gamma}}^+$ and $\psi_{{\gamma}}^-$ that :
\begin{equation}\label{efid}
 \psi_{{\gamma}}^+(i\lambda)\;\psi_{\gamma}^-(i\lambda)=\frac{\gamma}{\gamma+\psi(\lambda)}.\\
\end{equation}
\noindent On the other hand, it has been shown by Rogozin in \cite{ro} that for $\lambda,\mu\in\Cb$ with $\re{\lambda}\leq 0$, $\re{\mu}\geq 0$, we have :
\begin{equation}\label{erphi}
\psi_{{\gamma}}^+(\lambda)=E[\exp(\lambda S_{\ef_\gamma})]\;\;\textrm{ and}\;\; \psi_{\gamma}^-(\mu)=E[\exp(\mu I_{\ef_\gamma})].\\
\end{equation}
\noindent Hence, the representation of $\gamma{(\gamma+\psi(\lambda))}^{-1}$  in \eqref{efid} is an infinitely divisible factorization.\\

\end{rmq}

\begin{rmq}
\noindent We can find formula \eqref{eftk1} in \cite{Ber} (chap.VI, exercise 1, p.182) expressed as follows :
\begin{equation}\label{ec1}
\int\limits_0^{+\infty} e^{-\lambda a}E\left[\exp(-\beta T_{[a,+\infty[}-\theta X_{T_{[a,+\infty[}} )\right]da =\frac{\kappa(\beta,\lambda + \theta)-\kappa(\beta,\theta)}{\lambda\kappa(\beta,\lambda+\theta)}, 
\end{equation}
where $\kappa(\cdot,\cdot)$ is defined by :
\begin{equation}\label{defkapa}
\exp\left(-\kappa(\beta,\theta)\right)=E\left[\exp(-\beta\bar{\tau}_1-\theta S_{\bar{\tau}_1})\right], 
\end{equation}
\noindent where $\bar{\tau}$ denotes the right-continuous inverse of the local time at level $0$ of the reflected process $S-X$.\\

\noindent The relation between formulae \eqref{eftk1} and \eqref{ec1} is simple to establish. In fact, using the following formula (\cite{Ber}, Chap. VI, Corollary 10, p. 165) :
\begin{equation}\label{eqfc} 
\kappa(\beta,\theta)=\kappa(1,0)\exp\left(\int\limits_0^{+\infty} dt\int\limits_{[0,+\infty[} (e^{-t}-e^{-\beta t-\theta x})t^{-1}P(X_t\in dx)\right),
\end{equation}
\noindent we prove that for $q,\lambda>0$, we have :
$$\psi_q^+(-\lambda)=\frac{\kappa(q,0)}{\kappa(q,\lambda)}.$$
 \end{rmq}

\noindent In the case of a stable process with index $\alpha$, the following result will play a key role : 
\begin{prop}[\cite{Ber}, Chap.VIII, Proposition 2, p.219]\label{propber}
 We set $\rho=P(X_1>0)$. If $\rho\in ]0,1[$, then there exists $k>0$ such that :
\begin{equation}\label{ab}
 P(S_1\leq x)\sim k x^{\alpha\rho}\;\;\textrm{ as }\;\;x\rightarrow 0+. 
\end{equation}
\end{prop}
\section{About the law of $(\tx,\xtx)$}
\subsection{The general case}
\noindent The purpose of this paragraph is to express the joint Laplace transform of $(T_{[x,+\infty[},X_{T_{[x,+\infty[}})$ for $x>0$ :\\
$$E[\exp(-\lambda\tx -\mu\xtx)].$$
\noindent Formula $\eqref{eftk1}$ gives an expression of the Laplace transform of the quantity we are interested in. So the problem is 
reduced to invert this Laplace transform.\\

\noindent First, we note that for $\gamma>\mu>0$, we have :
$$\frac{1}{\gamma-\mu}=\int\limits_0^{+\infty}e^{-(\gamma-\mu) x}dx\quad\textrm{and}\quad\frac{\gamma}{\gamma-\mu}=1+\frac{\mu}{\gamma-\mu}.$$
\noindent On the other hand, from \eqref{erphi}, we obtain for $\gamma,\lambda>0$ :
\begin{equation}\label{erphi2}
\psi_{{\gamma}}^+(-\lambda)=E[\exp(-\lambda S_{\ef_\gamma})]=\int\limits_0^{+\infty}\lambda\; e^{-\lambda z}P(S_{\ef_\gamma}\leq z)\;dz.
\end{equation}
As a consequence, \eqref{eftk1} can be expressed as :
\begin{align}\label{eitx}
\izi e^{-\gamma x}\, E[&\exp(-\lambda\tx -\mu K_x)]\,dx\nonumber\\ &=\frac{1}{\gamma-\mu} -
\frac{1}{\psip{\lambda}{\mu}}\left(\frac{\gamma}{\gamma-\mu}\int\limits_0^{+\infty}e^{-\gamma z}P(S_{\ef_\lambda}\leq z)dz\right)\nonumber\\
&=\izi e^{-\gamma x}\left(e^{\mu x}-
\frac{P(\sel\leq x)}{\psip{\lambda}{\mu}}\right)dx-\frac{I(\gamma,\lambda,\mu)}{\psip{\lambda}{\mu}},
\end{align}
\noindent where :
$$I(\gamma,\lambda,\mu)=\mu\izi P(\sel \leq z)\left(\izi e^{-\gamma (y+z)}e^{\mu y} dy\right)dz.$$
\noindent Note that thanks to the change of variables $v=z$, $x=y+z$ and Fubini's theorem, we get :
$$I(\gamma,\lambda,\mu)= \mu \izi e^{-(\gamma-\mu) x}\left(\id{0}{x}e^{-\mu v}P(\sel\leq v) dv\right) dx.$$
\noindent If we put this expression for $I(\gamma,\lambda,\mu)$ in \eqref{eitx}, we obtain using \eqref{erphi2} :
\begin{equation}
\izi e^{-\gamma x} E[\exp(-\lambda\tx -\mu K_x)]dx=\frac{1}{\psip{\lambda}{\mu}}\izi e^{-(\gamma-\mu) x}J(\mu,\lambda;x)dx,
\end{equation}
\noindent where :
\begin{equation}\label{erj1}
J(\mu,\lambda;x)=\id{x}{+\infty}\mu e^{-\mu v}P(\sel\leq v) dv-e^{-\mu x}P(\sel\leq x).
\end{equation}
\noindent Therefore we have :
\begin{equation}\label{etxi1}
E[\exp(-\lambda\tx -\mu K_x)]=\frac{e^{\mu x}}{E\left[ e^{-\mu \sel}\right]}J(\mu,\lambda;x).
\end{equation}
\noindent Observe that :
\begin{align*}
 J(\mu,\lambda;x)&=E\left[\,\,\id{x\vee\sel}{+\infty}\mu e^{-\mu v}dv\right]-e^{-\mu x}P(\sel\leq x)\\
&=E\left[e^{-\mu (x\vee\sel)}\right]-e^{-\mu x}P(\sel\leq x),\\
\end{align*}
\noindent and thus :
\begin{equation}\label{erj2}
 J(\mu,\lambda;x)=E\left[ e^{-\mu \sel}1_{\{\sel\geq\, x\}}\right].
\end{equation}

\begin{thm}\label{tlctxxtx}
 \noindent For all $\lambda,\mu>0$ and $x\geq 0$ :
\begin{equation}\label{eqtxi2}
 E[\exp(-\lambda\tx -\mu\xtx)]
=\frac{E\left[ e^{-\mu \sel}1_{\{\sel\geq\, x\}}\right]}{E\left[ e^{-\mu \sel}\right]}.
\end{equation}
\end{thm}

\begin{proof}
It is immediate from \eqref{etxi1} and \eqref{erj2}.
\end{proof}

\subsection{The stable case}\label{ssstable}

\noindent Henceforth, $X$ will denote a real-valued stable Lévy process with index $\alpha\in]1,2]$. We denote : $$\rho=P(X_1\geq 0).$$
\noindent We know that in this case, $\rho\in[1-1/\alpha,1/\alpha]$ (see \cite{Ber}), in particular $\rho\in]0,1[$. The cases $\rho=1-1/\alpha$ and 
$\rho=1/\alpha$ correspond to the cases when $X$ has no negative jumps (spectrally positive case) and no positive jumps (spectrally negative case), respectively.

\subsubsection{Scaling for $\sel$ and some asymptotic results}
\noindent In Theorem \ref{tlctxxtx}, we see a link between the random variable $\sel$ and the joint distribution of $\left(\tx,\xtx\right)$. In this 
section, we show how the scaling property allows to study the absolute continuity of the law of $\sel$ and the asymptotic behavior of some quantities associated 
to this random variable.\\
\begin{prop}\label{ldsel}
For all $\lambda>0$, the law of the random variable $\sel$ is absolutely continuous. Its density $f_\lambda$ can be expressed as :
\begin{equation}
f_\lambda(x)=\frac{\lambda\alpha}{x}E\left[\tx \exp(-\lambda\tx)\right],\qquad x>0.
\end{equation}
\end{prop}
\begin{proof}
\noindent If we take the limit when $\mu$ goes to $0+$ in \eqref{eqtxi2}, we get :
$$P(\sel\leq x)=1-E\left[\exp(-\lambda\tx)\right].$$
\noindent On the other hand, the scaling property yields :
$$E\left[\exp(-\lambda\tx)\right]=E\left[\exp(-\lambda x^{\alpha}\tx[1])\right].$$
\noindent The result follows.
\end{proof}

\begin{lemme}[Scaling]\label{lsc}For all $\lambda,\mu>0$ and $x\geq 0$ :
\item (i)
$P(\sel\leq x)=P\left(\se\leq x\lambda^{1/\alpha}\right)=\izi e^{-v}P\left(S_1\leq\frac{x\lambda^{1/\alpha}}{v^{1/\alpha}}\right)dv.$
\item (ii) $E\left[ e^{-\mu \sel}\right]=E\left[ e^{-\frac{\mu}{\lambda^{1/\alpha}} \se}\right]=\izi e^{-v}P\left(\se\leq\frac{v\lambda^{1/\alpha}}{\mu}\right)dv.$
\end{lemme}
\begin{proof}
\item (i) We have for any $\lambda,\mu>0$ and $x\geq 0$. :
\begin{align*}
 P(\sel\leq x)&=\izi\lambda e^{-\lambda u}P(S_u\leq x)du\\
\textrm{(change of variables $\lambda u=v$)}&=\izi e^{-v}P(S_{v/\lambda}\leq x)dv\\
\textrm{(scaling)}&=\izi e^{-v}P\left(S_1\leq\frac{x\lambda^{1/\alpha}}{v^{1/\alpha}}\right)dv,
\end{align*}
\noindent which proves (i).\\ 
\item (ii) From \eqref{erphi2}, (i) and the change of variables $\mu z=v$ :
$$E\left[ e^{-\mu \sel}\right]=\izi\mu e^{-\mu z}P\left(\se\leq z\lambda^{1/\alpha}\right)dz
=\izi e^{-v}P\left(\se\leq\frac{v\lambda^{1/\alpha}}{\mu}\right)dv,$$

\noindent which completes the proof.
\end{proof}
\begin{prop}\label{altz}
There exists a constant $k^*>0$ such that :\\
\item(i) $P(\sel\leq x)\sim k^*\;x^{\alpha\rho}\;\lambda^\rho$ when $\lambda\rightarrow 0+$.\\
\item(ii) $\id{x}{+\infty}\mu e^{-\mu v}P(\sel\leq v) dv\sim k^*\left(\id{\mu x}{+\infty} e^{-y}y^{\alpha\rho} dy\right)\mu^{-\alpha\rho}{\lambda}^\rho$ when $\lambda\rightarrow 0+$.\\
\item(iii) $E\left[ e^{-\mu \sel}\right]\sim k^*\;\Gamma(1+\alpha\rho)\;\mu^{-\alpha\rho}\;{\lambda}^\rho$ when $\lambda\rightarrow 0+$.\\
\item(iv) $E\left[ e^{-\mu \sel}1_{\{\sel\geq\, x\}}\right]\sim k^*\alpha\rho\left(\id{\mu x}{+\infty} e^{-y}y^{\alpha\rho-1} dy\right)\mu^{-\alpha\rho}{\lambda}^\rho$ when $\lambda\rightarrow 0+$.
\end{prop}
\begin{proof}
\item (i) Thanks to Proposition \ref{propber} and Lemma \ref{lsc} (i), the dominated convergence theorem yields :
$$\frac{P(\sel\leq x)}{{\lambda}^\rho}=\izi e^{-v}\frac{P\left(S_1\leq\frac{x\lambda^{1/\alpha}}{v^{1/\alpha}}\right)}{{\lambda}^\rho}dv\underset{\lambda \rightarrow 0+}{\longrightarrow}k\;\Gamma(1-\rho)\;x^{\alpha\rho}.$$
\noindent It suffices to choose $k^*=k\,\Gamma(1-\rho)$ to finish the proof of (i). 
\item (ii) Applying the dominated convergence theorem, we deduce from (i) :
$$\frac{\id{x}{+\infty}\mu e^{-\mu v}P(\sel\leq v) dv}{{\lambda}^\rho}\underset{\lambda \rightarrow 0+}{\longrightarrow}k^*\id{x}{+\infty}\mu e^{-\mu v}v^{\alpha\rho} dv.$$
\noindent we obtain (ii) by the change of variables $y=\mu v$.
\item (iii) It suffices to consider $x=0$ in assertion (ii).
\item (iv) Using formulae \eqref{erj1} and \eqref{erj2}, (iv) is a consequence of (i) and (ii).
\end{proof}
\noindent \noindent In \cite{Ray}, Ray gives an expression for the density of the law of $\xtx$ in the symmetric case. In \cite{Bin}, Bingham generalizes this result 
to the case when $X$ is not spectrally negative. Now, we recover this result as a corollary of Theorem \ref{tlctxxtx} and Proposition \ref{altz} (see also Lemma 4.1 in \cite{yayayor}).\\
\begin{cor}[\cite{Ray}, \cite{Bin}]\label{cdxtx} If $\alpha\rho<1$, then :
\begin{equation}\label{lxtx}
\xtx    \overset{(loi)}{=}  \frac{x}{\beta_{{\alpha\rho},1-{\alpha\rho}}}.
\end{equation}
Consequently :
$$P(\xtx\in dy)=\rho(x,y)dy\textrm{, when } \rho(x,y)=\frac{\sin(\pi\alpha\rho)}{\pi}\frac{1}{y}{\left(\frac{x}{y-x}\right)}^{\alpha\rho}1_{]x,\infty[}(y).$$
\end{cor}
\begin{proof}
 \noindent If we let $\lambda$ tend to $0+$ in \eqref{eqtxi2}, we get with the help of Proposition \ref{altz},
\begin{equation*}
 E[\exp(-\mu\xtx)]=\frac{1}{\Gamma({\alpha\rho})} \id{\mu x}{+\infty} e^{-y}y^{\alpha\rho-1} dy,
\end{equation*}
\noindent and by using the identity $y^{\alpha\rho-1}=\frac{1}{\Gamma(1-\alpha\rho)}\izi e^{-yv}v^{-\alpha\rho}dv$, we obtain :
\begin{align*}
 E[\exp(-\mu\xtx)]&=\frac{\sin(\pi\alpha\rho)}{\pi} \id{\mu x}{+\infty} dy\izi e^{-y(v+1)}v^{-\alpha\rho}dv\\
&=\frac{\sin(\pi\alpha\rho)}{\pi}\izi \frac{e^{-\mu x (v+1)}}{v+1}v^{-\alpha\rho}dv\quad\textrm{(Fubini's theorem)}\\
&=\frac{\sin(\pi\alpha\rho)}{\pi}\id{0}{1}e^{-\mu x/z}z^{\alpha\rho -1}{(1-z)}^{-\alpha\rho}dz\quad\textrm{(change of}\\ 
&\qquad\qquad\qquad\qquad\qquad\quad\quad\textrm{ variables } z=1/(v+1)),
\end{align*}

\noindent which proves Corollary \ref{cdxtx}.
\end{proof}
\begin{rmq}
\noindent Similarly, the law of $K_x$ is absolutely continuous. Its density is given by :
$$\rho_{K}(x,y)=\rho(x,x+y)=\frac{\sin(\pi\alpha\rho)}{\pi}\frac{1}{(y+x)}{\left(\frac{x}{y}\right)}^{\alpha\rho}1_{]0,\infty[}(y).$$
\end{rmq}
\begin{rmq}
\noindent This result has been generalized in the papers \cite{Doky} and \cite{Kypar}.
\end{rmq}

\subsubsection{An asymptotic law}
\noindent Note first that a standard application of the Markov property permits to show that :
$$ P(S_t\in\,]\,x,x+h\,])=\!\!\!\!\!\!\!\iint\limits_{[0,t[\times ]x,x+h]}\!\!\!\!\!\!\! P(T_{[x+h-y,+\infty[}\geq t-s)\,P(\tx\in ds, \xtx\!\in dy).$$
\noindent Taking the limit in this equality when $h$ tends to $0$, we may expect to find links between the absolute continuity of the law of $S_t$ and the asymptotic behavior of the 
random variable $\tx$ conditioned to $\{K_x\leq h\}$.\\

\noindent With this motivation, we introduce for each $x,h>0$ the random variable $\txh$ whose law is given by :\\
$$P(\txh\in\cdot\,\,)=P(\tx\in\cdot\,\,\,|\,K_x\leq h).$$

\noindent Thus the aim in this section is to study the asymptotic behavior of the variables $\txh$ when $h$ goes to $0+$. For this, we start by 
computing the asymptotic probability of the  events $\{K_x\leq h\}$.\\
\begin{lemme}\label{ckxp} For all $x>0$ :
\begin{equation}
\frac{P(K_x\leq h)}{h^{1-\alpha\rho}}\underset{h \rightarrow 0+}{\longrightarrow} \frac{\sin (\pi\alpha\rho)}{\pi(1-\alpha\rho)}\;x^{\alpha\rho -1}.
\end{equation}
\end{lemme}
\begin{proof}
\noindent By Corollary \ref{cdxtx}, we have :
\begin{align*}
 P(\xtx\leq x+h) &= \frac{\sin(\pi\alpha\rho)}{\pi}\id{x}{x+h}\frac{1}{y}{\left(\frac{x}{y-x}\right)}^{\alpha\rho}dy\\
&=\frac{\sin(\pi\alpha\rho)}{\pi}\;\;x^{\alpha\rho}h^{1-\alpha\rho}\id{0}{1}\frac{1}{(uh+x)u^{\alpha\rho}}du\quad\textrm{(change of}\\
&\qquad\qquad\qquad\qquad\qquad\qquad\textrm{ variables }u=(y-x)/h).
\end{align*}
\noindent The result follows by dividing both sides by $h^{1-\alpha\rho}$ and applying the dominated convergence theorem.
\end{proof}
\noindent Now, we define the probability measure $\plx$ by : $$\plx(A)=\frac{E\left[1_A \exp(-\lambda\tx)\right]}{E\left[\exp(-\lambda\tx)\right]},\quad A\in\mathcal{F}_\infty.$$
\noindent The following lemma is the analogue of Lemma \ref{ckxp} for the probability measure $\plx$.
\begin{lemme}\label{ckxpplx} For all $\lambda,x>0$, we have :
\begin{equation}
\frac{\plx(\xtx-x\leq h)}{h^{1-\alpha\rho}}\underset{h \rightarrow 0+}{\longrightarrow} \frac{\sin (\pi\alpha\rho)}{k^*\pi\alpha\rho(1-\alpha\rho)}\;\frac{\lambda^{-\rho}f_\lambda (x)}{P(\sel\geq x)},
\end{equation}
where $k^*$ is the constant appearing in Proposition \ref{altz}.
\end{lemme}
\begin{proof}
We consider the fonction $U:[0,+\infty[\rightarrow[0,+\infty[$ defined by : 
$$U(h)=\plx\left(K_x\leq h\right).$$ 
\noindent By Tauberian theorem (see p.10 in \cite{Ber}), the behavior of $U$ around $0+$ is related to the behavior of its Laplace transform at infinity.\\

\noindent We remark that :
\begin{equation}\label{lpu}
\izi e^{-\mu y}U(dy)=\elx\left[\exp(-\mu K_x)\right]=\frac{E\left[\exp(-\lambda\tx-\mu K_x)\right]}{E\left[\exp(-\lambda\tx)\right]}.
\end{equation}
On the other hand, thanks to \eqref{erj1}, \eqref{etxi1} and an obvious change of variables, we obtain :
\begin{equation}\label{lpu1}
E\left[\exp(-\lambda\tx-\mu K_x)\right]=\frac{\izi e^{-y}\left(P\left(\sel\leq \frac{y}{\mu}+x\right)-P\left(\sel\leq x\right)\right)dy}{
E\left[e^{-\mu\sel}\right]}.
\end{equation}
Using Lemma \ref{ldsel} and applying the dominated convergence theorem, we show that :
\begin{equation}\label{am1}
 \mu\izi e^{-y}\left(P\left(\sel\leq \frac{y}{\mu}+x\right)-P\left(\sel\leq x\right)\right)dy\underset{\mu \rightarrow +\infty}{\longrightarrow}f_\lambda(x).
\end{equation}
We know from part (ii) of Lemma \ref{lsc} that :
$$E\left[e^{-\mu\sel}\right]=E\left[e^{-\frac{1}{\lambda^{1/\alpha}}\sel[\mu^{-\alpha}]}\right],$$
and therefore, from part (iii) of Proposition \ref{altz}, we obtain that :
\begin{equation}\label{am2}
 \mu^{\alpha\rho}E\left[e^{-\mu\sel}\right]\underset{\mu \rightarrow +\infty}{\longrightarrow}k^*\Gamma(1+\alpha\rho)\lambda^\rho.
\end{equation}
Thus, from \eqref{lpu}, \eqref{lpu1}, \eqref{am1} and \eqref{am2}, we get :
\begin{equation}
 \mu^{1-\alpha\rho}\izi e^{-\mu y}U(dy)\underset{\mu \rightarrow +\infty}{\longrightarrow}\frac{1}{k^*\Gamma(1+\alpha\rho)}\frac{\lambda^{-\rho}f_\lambda (x)}{P(\sel\geq x)},
\end{equation}
and then, thanks to the Tauberian theorem (see p.10 in \cite{Ber}), we obtain :
\begin{equation}
\frac{1}{h^{1-\alpha\rho}}\;U(h)\underset{h \rightarrow 0+}{\longrightarrow}\frac{1}{k^*\Gamma(1+\alpha\rho)\Gamma(2-\alpha\rho)}\frac{\lambda^{-\rho}f_\lambda (x)}{P(\sel\geq x)},
\end{equation}
which completes the proof.
\end{proof}
\begin{prop}\label{cdtl}
For every $\lambda,x>0$, we have :
\begin{equation}
E\left[\exp(-\lambda\tx\;)\;|\;K_x\leq h\right]\underset{h \rightarrow 0+}{\longrightarrow}\frac{1}{k^*\alpha\rho}\,x^{1-\alpha\rho}\,\lambda^{-\rho}\,f_\lambda(x).
\end{equation}
\end{prop}
\begin{proof}
We remark that :
\begin{equation*}
 E\left[\exp(-\lambda\tx\;)\;|\;K_x\leq h\right]=\frac{\plx\left(K_x\leq h\right)}{P\left(K_x\leq h\right)}E\left[\exp(-\lambda\tx)\right].
\end{equation*}
\noindent Now, the result is a consequence of Lemmas \ref{ckxp} and \ref{ckxpplx}.
\end{proof}
\begin{lemme}\label{ppm}
For all $x>0$ :
\begin{equation}
\lim\limits_{\lambda\rightarrow 0+} \frac{1}{k^*\alpha\rho}\,x^{1-\alpha\rho}\,\lambda^{-\rho}\,f_\lambda(x)=1.
\end{equation}
\end{lemme}
\begin{proof}
\noindent The scaling property entails :
\begin{equation*}
E[\tx\exp(-\lambda\tx)]=\izi e^{-v}(1-v)P\left(S_1\leq\frac{x\lambda^{1/\alpha}}{v^{1/\alpha}}\right)\,dv, 
\end{equation*}
\noindent and then, by definition of $f_\lambda$ we obtain from Proposition \ref{propber} and the dominated convergence theorem that :
\begin{align*}
\lambda^{-\rho}\,f_\lambda(x)& \underset{\lambda \rightarrow 0+}{\longrightarrow} k\,\alpha \,x^{\alpha\rho-1}\,\izi e^{-v}\,(1-v)\,v^{-\rho}\, dv\\
&=k\,\alpha\,\rho\,\Gamma(1-\rho)\,x^{\alpha\rho-1}.
\end{align*}
\noindent The result follows.
\end{proof}
\begin{thm}\label{thmll}
For each $x>0$, the family of random variables ${\{\txh\}}_{h>0}$ converges in law as $h$ tends to $0+$. The limit is denoted by $\txz$ and its law is given by :
$$P\left(\,\txz\leq\, t\,\right)=\frac{\sin (\pi\rho)}{k\,\pi\rho}\,x^{-\alpha\rho}E\left[\frac{\tx \,\,1_{\{\tx\leq\, t\}}}{{\left(t-\tx\right)}^{1-\rho}}\right].$$ 
\end{thm}
\begin{proof}
\noindent For $x,h>0$, we denote by $\Ls_h^x$ the Laplace transform of the random variable $\txh$ and by $\Ls^x$ the fonction defined by :
$$\Ls^x(\lambda)=\frac{1}{k^*\alpha\rho}\,x^{1-\alpha\rho}\,\lambda^{-\rho}\,f_\lambda(x).$$   

\noindent According to Proposition \ref{cdtl}, the Laplace transform $\Ls_h^x$ converges pointwise to the fonction $\Ls^x$. Since we have already 
shown in Lemma \ref{ppm} that :
$$\lim\limits_{\lambda\rightarrow 0+}\Ls^x(\lambda)=1,$$
it appears that the fonction $\Ls^x$ is the Laplace transform to some probability measure with support in $\Rb_+$ (see \cite{Chung}, Theorem 6.6.3, p.190). 
We denote this limit measure by $\nu_x$.\\

\noindent We remark that :
\begin{equation}\label{els1}
 \frac{\Ls^x(\lambda)}{\lambda}= \frac{1}{\lambda}\int\limits_{[0,+\infty]}e^{-\lambda y}\nu_x(dy)=\izi e^{-\lambda y}\nu_x([0,y])\,dy.
\end{equation}
\noindent On the other hand, we have :
$$\frac{\Ls^x(\lambda)}{\lambda}=\frac{1}{k^*\rho}\,x^{-\alpha\rho}\,\lambda^{-\rho}E\left[\tx \exp(-\lambda\tx)\right].$$
By using the identity $\lambda^{-\rho}=\frac{1}{\Gamma(\rho)}\izi e^{-\lambda z}z^{\rho-1}dz$, we get :
\begin{equation}\label{els2}
\frac{\Ls^x(\lambda)}{\lambda}=\frac{\sin (\pi\rho)}{k\,\pi\rho}\,x^{-\alpha\rho}\izi e^{-\lambda u}E\left[\frac{\tx \,\,1_{\{\tx\leq\, u\}}}{{\left(u-\tx\right)}^{1-\rho}}\right]\, du.
\end{equation}
The result is obtained by comparing \eqref{els1} and \eqref{els2}.
\end{proof}
\subsubsection{The laws of $\tx$ and $\txz$ have no point masses}
\noindent In  \cite{pero} (Lemma 1), it has been demonstrated in a more general framework that the law of $S_t$ admits no point masses, i.e., for any $x\geq 0$, $P(S_t=x)=0$. Since we have $\{S_t\geq x\}=\{\tx\leq t\}\cup\{S_t=x\}$, then :
$$P(S_t\geq x)=P(\tx\leq t).$$
So, by using the scaling property, we obtain in particular :
$$S_1^{-\alpha}\overset{(loi)}{=}\tx[1].$$
Thus we see that $\tx[1]$ admits no point masses either (the same for $\tx$, by scaling). In the following proposition, we find this result directly from the scaling property.
\begin{prop}\label{pmp}
For every $x,t>0$ :
$$P(\tx=t)=P(\txz=t)=0.$$
\end{prop}

\begin{proof}
\noindent We first show the result for $\tx$. By scaling property, with no loss of generality, we may suppose that $x=1$.\\

\noindent Suppose, by contradiction, that there is $t_0>0$ such that :
$$\delta:=P(\tx[1]=t_0)>0.$$ 
\noindent Suppose now that there is $h_0> 0$ satisfying :
$$\eta:=P(\tx[1]=t_0,\, K_1>h_0)>0.$$
\noindent In this case, we have :
$$P(\forall u\in[1,1+h_0/2] \,:\,\tx[u]=t_0)>\eta,$$
\noindent and therefore, for every $u\in[1,1+h_0/2]$, we have $P(\tx[u]=t_0)>\eta$. Thus, the scaling property yields that for each $u\in[1,1+h_0/2]$ : $$P(\tx[1]=t_0/u^\alpha)>\eta.$$ 
\noindent In other words, for all $s\in [(2/(2+h_0))^\alpha t_0,t_0]$ : $$P(\tx[1]=s)>\eta,$$
\noindent which can not be true. We have therefore, for any $h>0$ :  \\
$$P(\tx[1]=t_0,\, K_1>h)=0.$$
\noindent  Thus, for every $h>0$ :
$$0<\delta=P(\tx[1]=t_0)=P(\tx[1]=t_0,\, K_1\leq h)\leq P(K_1\leq h),$$
\noindent which is a contradiction, because the right-hand side converges to $0$ when $h$ tends to $0$.\\

\noindent The result for $\txz$ can be obtained from the result for $\tx$, since thanks to Theorem \ref{thmll} we have :
 
\begin{equation*}
 P(t-h<\txz\leq t)\leq\frac{\sin (\pi\rho)}{k\,\pi\rho}\,x^{-\alpha\rho}E\left[\frac{\tx \,\,1_{\{t-h<\tx\leq\, t\}}}{{\left(t-\tx\right)}^{1-\rho}}\right].\\
\end{equation*}
\end{proof}
\noindent\textbf{Comment.} In a forthcoming publication, I would like to establish more precise links between the absolute continuity of the law of $S_t$ and the asymptotic random variable $\txz$.
\bibliographystyle{plain}
\bibliography{Reference}

\end{document}